\documentclass[11pt]{article}
\usepackage{amscd, amsmath, amssymb}

\title{Moving codimension-one subvarieties over finite fields}
\author{Burt Totaro}
\date{  }

\def\Z{\text{\bf Z}}

\def\R{\text{\bf R}}

\def\P{\text{\bf P}}
\def\F{\text{\bf F}}
\def\Fb{\overline{\F}}
\def\arrow{\rightarrow}
\def\inj{\hookrightarrow}

\def\qed{\ QED }

\def\ker{\text{ker}\;}
\def\Pic{\text{Pic}}

\def\Ext{\text{Ext}}

\begin{document}
\maketitle

\newtheorem{theorem}{Theorem}[section]
\newtheorem{corollary}[theorem]{Corollary}
\newtheorem{lemma}[theorem]{Lemma}

In topology, the normal bundle of a submanifold
determines a neighborhood of the submanifold up to isomorphism.
In particular, the normal bundle of a codimension-one
submanifold is trivial if and only if the submanifold can be moved
in a family of disjoint submanifolds. In algebraic geometry,
however, there are higher-order obstructions to moving
a given subvariety.

In this paper, we develop an obstruction theory, in the spirit
of homotopy theory, which gives some
control over when a codimension-one subvariety moves in a family of disjoint
subvarieties. Even if a subvariety
does not move in a family, some positive multiple of it may.
We find a pattern linking the infinitely many obstructions to moving
higher and higher multiples of a given subvariety. As an application,
we find the first examples of line bundles $L$ on smooth projective
varieties over finite fields which are nef ($L$ has nonnegative
degree on every curve) but not semi-ample (no positive power of $L$
is spanned by its global sections). This answers
questions by Keel and Mumford.

Determining which line bundles are spanned by their global
sections, or more generally are semi-ample, is a fundamental
issue in algebraic geometry. 
If a line bundle $L$ is semi-ample, then the powers of $L$ determine
a morphism from the given variety onto some projective variety.
One of the main problems of the minimal model program,
the abundance conjecture, predicts that a variety with
nef canonical bundle should have semi-ample canonical bundle
\cite[Conjecture 3.12]{KM}.

One can hope to get more insight into the abundance conjecture
by reducing varieties in characteristic zero
to varieties over finite fields,
where they become simpler in some ways.
In particular, by
Artin, every nef line bundle $L$ with $L^2>0$ on a projective surface
over the algebraic closure
of a finite field is semi-ample
\cite[proof of Theorem 2.9(B)]{Artin}.
This is far from true over other algebraically closed fields,
by Zariski \cite[section 2]{Zariski}.
Keel generalized Artin's theorem, giving powerful sufficient conditions
for a nef line bundle on a projective variety over $\Fb_p$
to be semi-ample \cite{KeelAnn,KeelPol}. As an application,
he constructed contractions of the moduli space of stable curves which
exist as projective varieties in every finite characteristic, but
not in characteristic zero.

Keel asked whether a nef line bundle $L$ on a smooth projective
surface over $\Fb_p$ is always semi-ample, the open case being line bundles
with $L^2=0$. (This is part of his Question 0.8.2 \cite{KeelPol},
in view of Theorem \ref{Leq} below.)

Using our obstruction theory, we can see where counterexamples
to Keel's Question 0.8.2 should be expected and produce them.
We obtain the first known examples
of nef but not semi-ample line bundles on smooth projective varieties
over $\Fb_p$, for any prime number $p$ (Theorem \ref{mainex}).
Equivalently, we give
faces of the closed cone of curves which have
rational slope but cannot be contracted.
The line bundles we construct are effective, of the form
$O(C)$ for a smooth curve $C$ of genus 2 with self-intersection zero
on a smooth projective surface.
Thus we give the first examples over $\Fb_p$ of a curve 
with self-intersection zero such that
no multiple of the curve moves in a family of disjoint curves
(which would give a fibration of the surface over a curve).
This answers a question raised by Mumford \cite[p.~336]{MumfordI}.

One can still hope for some positive
results in the direction of Question 0.8.2. Sakai \cite{SakaiComp}
and Ma\c sek \cite{Masek} gave a positive result when the curve
$C$ has genus at most 1; we give Keel's proof of their result in
Theorem \ref{genus1}.
(This makes sense in terms of minimal model theory,
since a curve $C$ with $C^2=0$ in a surface $X$ has genus
at most 1 exactly when $K_X\cdot C\leq 0$.)

Building upon our basic example, we exhibit a nef {\it and big }line
bundle on a smooth projective 3-fold over $\Fb_p$ which is not
semi-ample (Theorem \ref{bigex}).

Keel pointed out the following special case of Question 0.8.2, which
remains open and very interesting \cite[Question 0.9]{KeelPol}.
Let $X$ be a smooth projective surface over $\Fb_p$ with a line bundle $L$.
If $L\cdot C>0$ for every curve $C$, does it follow that $L$ is ample
(or equivalently, that $L^2>0$)? Counterexamples to this statement,
using ruled surfaces,
were given over the complex numbers by Mumford
\cite[Example 10.6]{HartshorneAmple},
and over uncountable algebraically closed fields
of positive characteristic by Mehta and Subramanian
\cite[Remark 3.2]{MehtaS}.

Thanks to Daniel Huybrechts,
Yujiro Kawamata, Sean Keel, and James McKernan for their comments.

\section{Notation}
\label{notation}

Varieties are reduced and irreducible by definition.
A {\it curve }on a variety means a
closed subvariety of dimension 1.
A line bundle $L$ on an $n$-dimensional proper variety $X$ over a field
is {\it nef }if the intersection number $L\cdot C$ is nonnegative 
for every curve $C$ in $X$. We often use additive notation
for line bundles, in which the line bundle $L^{\otimes a}$
is called $aL$ for an integer $a$. For a nef line bundle $L$, a curve $C$
is called {\it $L$-exceptional }if $L\cdot C=0$.
A line bundle $L$ is {\it big }if the rational map associated
to some positive multiple of $L$ is birational. We use
the following fact \cite[Theorem VI.2.15]{KollarRat}.

\begin{lemma}
\label{big}
Let $L$ be a nef line bundle on a proper variety of dimension $n$
over a field.
Then $L$ is big if and only if $L^n>0$.
\end{lemma}

A line bundle $L$ is {\it semi-ample }if
some positive multiple of $L$ is spanned by its global sections.
A semi-ample line bundle is nef. Also, a semi-ample line
bundle determines a contraction: the algebra $R(X,L)=
\oplus_{a\geq 0} H^0(X,aL)$
is finitely generated by Zariski \cite{Zariski}, $Y:=\text{Proj }R(X,L)$
is a projective variety, and there is a natural surjection $f$
from $X$ to $Y$ with connected fibers (meaning that $f_*O_X=O_Y$).
Conversely, any surjective morphism with connected fibers
from $X$ to a projective variety $Y$ arises
in this way from some semi-ample line bundle $L$ on $X$, by taking $L$
to be the pullback of any ample line bundle on $Y$.

Mourougane and Russo gave a useful decomposition of semi-ampleness
into two properties,
extending earlier results by Zariski \cite{Zariski}
and Kawamata \cite{Kawamata}. Let $L$ be a nef
line bundle on a variety $X$.
Define the {\it
numerical dimension }$\nu(X,L)$ to be the largest
natural number $\nu$ such that the cycle $L^{\nu}$ is numerically
nontrivial, that is, the largest $\nu$ such that $L^{\nu}\cdot S>0$
for some subvariety $S$ of dimension $\nu$.
The {\it Iitaka dimension }$\kappa(X,L)$ is defined to be
$-\infty$ if $H^0(X,aL)=0$ for all $a>0$. Otherwise, we define $\kappa(X,L)$
to be one less than
the transcendence degree of the quotient field of the graded algebra $R(X,L)$,
or equivalently the maximum dimension of the image of $X$ under
the rational maps to projective space defined by powers of $L$.
We always have
$\kappa(X,L)\leq \nu(X,L)$, and $L$ is called {\it good }if the two
dimensions are equal.

\begin{theorem} \cite[Corollary 1]{MouR}
\label{good}
Let $X$ be a normal proper variety over a field $k$.
A nef line bundle $L$ on $X$ is semi-ample if and only if it is good
and the graded algebra $R(X,L)$ is finitely generated over $k$.
The finite generation is automatic when $\kappa(X,L)\leq 1$.
\end{theorem}

For example, when $L^n>0$, in other words when $L$ has maximal numerical
dimension, Lemma \ref{big} says that $L$ is automatically good,
and so semi-ampleness is purely a question of finite generation.
This paper is about the ``opposite'' situation, where $L$
has numerical dimension 1; by Theorem \ref{good},
semi-ampleness is equivalent to goodness in this case. On the other hand,
problems of semi-ampleness always have some relation to problems of
finite generation. For example, if $L$ is a nef line bundle on
a projective variety $X$ and $M$ is an ample line bundle, then
$L$ is semi-ample if and only if the algebra
$\oplus_{a,b\geq 0}H^0(X,aL+bM)$ is finitely generated, as one easily
checks.

Finally, here is one last standard result.

\begin{lemma}
\label{torsion}
Every numerically
trivial line bundle on a proper scheme over the algebraic closure
of a finite field is torsion.
\end{lemma}

The proof is based on the fact that an abelian variety
has only finitely many
rational points over any given finite field. A reference is 
\cite[Lemma 2.16]{KeelAnn}.

\section{Moving elliptic curves on surfaces over finite fields}

Sakai \cite[Theorem 1, Proposition 5,
and Concluding Remark]{SakaiComp}
and Ma\c sek \cite[Theorem 1 and Lemma]{Masek} showed that
for any curve of genus 1
with self-intersection zero
on a smooth projective surface over $\Fb_p$, some positive multiple
of the curve always moves. In this section we give
Keel's proof of this result.
The corresponding statement is easy for curves
of genus 0 on a surface (over any field, in fact)
and false for genus at least 2 (Theorem \ref{mainex}).

\begin{theorem}
\label{genus1}
Let $C$ be a curve of arithmetic genus 1 in a smooth projective surface
$X$ over $\Fb_p$. If $C^2=0$, then $L=O(C)$ is semi-ample.
Equivalently, $C$ is a fiber (possibly a multiple fiber)
in some elliptic or quasi-elliptic fibration of $X$.
\end{theorem}

This is false over all algebraically closed fields $k$ except the algebraic
closure of a finite field, by some standard examples.
First, let $C$ be a smooth cubic curve in
$\P^2$, and let $X$ be the blow-up of $\P^2$ at 9 points on the curve.
Since $C^2=9$ in $\P^2$, the proper transform $C$ in $X$ has $C^2=0$.
By appropriate choice of the points to blow up, we can make the normal
bundle of $C$ in $X$ any line bundle of degree 0 on $C$. If $k$
is not the algebraic closure of a finite field, this normal bundle
can be non-torsion in the Picard group of $C$,
and so no positive multiple of $C$
moves in $X$. Another example, on a ruled surface, works
only in characteristic zero: see the discussion after Lemma
\ref{mumford}.

{\bf Proof. }Since the normal bundle $N_{C/X}=L|_C$ has degree 0,
it is torsion by Lemma \ref{torsion} since the base field $k$
is $\Fb_p$. Let $m$ be the order of $L|_C$ in $\Pic(C)$. For any
integer $a$, we have the exact sequence:
$$H^0(C,aL)\arrow H^1(X,(a-1)L)\arrow H^1(X,aL)\arrow H^1(C,aL).$$
Since $C$ has genus 1 and the line bundle $aL$ is nontrivial
on $C$ for $1\leq a\leq m-1$, we have $H^0(C,aL)=0$ and 
$H^1(C,aL)=0$ in that range. So the exact sequence gives isomorphisms
$$H^1(X,O)\cong H^1(X,L)\cong \cdots \cong H^1(X,(m-1)L).$$

Next, we have the exact sequence
$$0\arrow H^0(X,(m-1)L)\arrow H^0(X,mL)\arrow H^0(C,mL)
\arrow H^1(X,(m-1)L).$$
Here $H^0(C,mL)\cong k$. Thus, if $H^1(X,O)$ is zero, then $H^1(X,
(m-1)L)$ is zero, and so the sequence shows that the divisor $mC$
moves nontrivially in $X$. Therefore $\kappa(X,L)\geq 1$,
and $L$ is semi-ample by Theorem \ref{good}.

In general, we use the idea of ``killing cohomology'':
for any projective variety $X$ over
a field $k$ of characteristic $p>0$, and any element $\alpha$
of $H^1(X,O)$, there is a surjective morphism $f:W\arrow X$ of projective
varieties which kills $\alpha$, meaning that $f^*(\alpha)=0$.
Indeed, by Serre, for $X$ smooth
one can do this with a finite flat morphism $f$
(a composite of etale $\Z/p$-coverings and Frobenius morphisms)
\cite[Proposition 12 and section 9]{SerreTop}. The construction shows
that $W$ is smooth since $X$ is. Then the above proof applied
to $W$ shows that $\kappa(W,f^*(L))\geq 1$. By Ueno, for any
surjective morphism $f$ of normal projective varieties, we have
$\kappa(W,f^*(L))=\kappa(X,L)$ \cite[Theorem 5.13]{Ueno}.
Thus $\kappa(X,L)\geq 1$
and so $L$ is semi-ample. \qed

\section{$L$-equivalence}

In this section, we show that $L$-equivalence (as defined by Keel)
is automatically
bounded on a normal projective surface over any field. The definition
is given before Theorem \ref{Leq}.
(By contrast, Koll\'ar observed
that $L$-equivalence can be unbounded
on a non-normal surface \cite[section 5]{KeelPol}.) 
The proof is elementary geometry of surfaces (the Hodge index
theorem).

The ``reduction map for nef line bundles'' of Bauer et
al.\ \cite{BC} is a similar application of the Hodge index theorem which
works in any dimension, but some work would be needed to go from
their theorem to the boundedness of $L$-equivalence on normal surfaces.
Their theorem is stated over the complex numbers, but the argument
works in any characteristic.
The difficulty is that their theorem only applies
to ``general'' points, meaning
points outside a countable union of proper subvarieties \cite[Theorem 2.1,
2.4.2]{BC}. For example, it seems to be unknown
whether there is a nef line bundle $L$
on some normal complex projective variety $X$ such that the set of curves
$C$ with $L.C=0$ is countably infinite. That does not
happen for $X$ of dimension 2 by Lemma \ref{bound} below
and \cite[Theorem 2.1]{Totarodivisor},
but the latter result uses the Albanese map in a way
that has no obvious generalization to higher dimensions.

Let $\rho(X)$ denote the Picard number of a projective variety $X$,
that is, the dimension of the real vector space $N^1(X)$ generated
by line
bundles modulo numerical equivalence.

\begin{lemma}
\label{bound}
Let $X$ be a smooth projective surface over a field. Let $L$
be a nef line bundle on $X$ such that $L^2=0$ and $L$ is numerically
nontrivial. Then there are at most $2(\rho(X)-2)$ curves $A$ on $X$
such that $L\cdot A=0$ and $A\not\in \R^{>0}\cdot L\subset N^1(X)$.
(These curves $A$ will all have negative self-intersection.)

If in addition there is no effective 1-cycle in $\R^{>0}\cdot L$,
or if there is a curve $C$ in $X$ such that every effective
1-cycle in $\R^{>0}\cdot L$ is a multiple of $C$ as a cycle,
then there are at most $\rho(X)-2$ curves $A\neq C$ with
$L\cdot A=0$.

Finally, let $L$ be a nef line bundle on a surface
with $L^2>0$. Then there are at most $\rho(X)-1$ curves $A$ with $L\cdot
A=0$. 
\end{lemma}

All these bounds are optimal, as we now check.

{\bf Example. }Let $X$ be the blow-up of $\P^1\times \P^1$
at $d$ points whose projections to the second factor are distinct.
Then $X$ has Picard number $d+2$.
Let $L$ be the pullback to $X$ of the line
bundle $O(1)$ on the second factor of $\P^1$. Then there are
exactly $2d$ curves $A$ with $L\cdot A=0$ which are not in
the ray $\R^{>0}\cdot L\subset N^1(X)$, namely the proper transforms
of the fibers containing the $d$ given points, together with the
$d$ exceptional curves. These $2d$ curves are all $(-1)$-curves.
This shows the optimality of the first statement in Lemma \ref{bound}.

Next, assume that the base field has characteristic zero, and
let $X$ be the ruled surface over a curve of genus at least 1
associated to a nontrivial extension of the trivial line bundle by
itself.
Let $C$ be the section in $X$ with zero self-intersection.
Then every curve in the ray $\R^{>0}\cdot C\subset N^1(X)$ is
a multiple of $C$ as a cycle.
In this case, $X$ has Picard number 2 and,
by Lemma \ref{bound}, there is no curve $A\neq C$ with
$C\cdot A=0$. Blowing up $X$ at $d$ points not on $C$
gives a surface $M$ with Picard number $d+2$ such that the line bundle
$L=O(C)$ on $M$ has exactly $d$ curves $A\neq C$ with $C\cdot A=0$,
namely the $d$ exceptional curves. This shows the optimality
of the second statement in Lemma \ref{bound}.

Finally, let $X$ be the blow-up of $\P ^2$ at $d$ points, and
let $L$ be the pullback of the line bundle $O(1)$ to $X$, which is nef and big.
Then $X$ has Picard number
$d+1$, and there are exactly $d$ curves on which $L$ has degree zero,
namely the $d$ exceptional curves. This shows the optimality of
the last statement of the lemma.

{\bf Proof of Lemma \ref{bound}. }We have that $L^2=0$ and
$L$ is numerically nontrivial. By the Hodge index theorem
\cite[Theorem V.1.9]{HartshorneAG},
the intersection form on $N^1(X)$ passes to a negative
definite form on $V:=L^{\perp}/(\R\cdot L)$. As a result,
every curve $A$ in $X$ with $L\cdot A =0$ and $A\not\in \R^{>0}\cdot
L$ has $A^2<0$. Also, for any two distinct curves $A_1$ and $A_2$
with these properties, we have $A_1\cdot A_2\geq 0$.

Let us change the sign of the intersection form on $V$, to make
it positive definite. Then the set $S$ of curves as above maps
injectively
to a subset of $V$ such that any two distinct elements of $S$
have inner product $\langle x,y\rangle \leq 0$. Rankin showed that
a subset $S$ of a real inner product space $V$ with these
properties has at most $2\dim(V)$ elements \cite{Rankin}.
The proof is quite short. In our case, $V$ has dimension
$\rho(X)-2$, and so $S$ has order at most
$2(\rho(X)-2)$, as we want.

We now prove the second statement of the lemma. Again,
change the sign of the intersection form on $V$ to make
it positive definite. In this case,
we have the additional information that $0\in V$
is not in the convex hull of the set $S$. So there is a vector
$u$ in $V$ with $\langle x,u \rangle >0$ for all $x\in S$. 
Adding $-\epsilon u$ to each $x\in S$ for a small positive number
$\epsilon $ changes the set $S$ so as to have
$\langle x,y\rangle <0$ (not just $\leq 0$) for all 
$x\neq y$ in $S$, while still having $\langle x,u\rangle >0$.
Then the set $S\cup \{-u\}$ has the property that the inner
product between any two of its elements is negative. In that case,
Rankin showed that $S\cup \{-u\}$ has at most $\dim(V)+1$
elements \cite{Rankin}. So $S$ has at most $\dim(V)=\rho(X)-2$
elements, as we want.

For the third statement of the lemma, we have $L^2>0$. In this case,
the Hodge index theorem gives that the intersection form on 
$L^{\perp}$ is negative definite. Let $V$ denote $L^{\perp}$
in this case, with the negative of that intersection form.
Here $0\in V$ is not
in the convex hull of the set $S$ of classes of curves $A$
with $L\cdot A=0$, using that $X$ is projective. Then the previous
paragraph's argument shows that $S$ has at most $\dim(V)=\rho(X)-1$
elements. \qed

We now recall Keel's notion of $L$-equivalence.
Let $L$ be a nef line bundle
on a proper scheme $X$ over a field. Two closed points
in $X$ are called {\it $L$-equivalent }if they can be connected
by a chain of curves $C$ such that $L\cdot C=0$.
Say that $L$-equivalence is {\it bounded }if there is a positive
integer $m$ such that any two $L$-equivalent points can be connected
by a chain with length at most $m$ of such curves.

Keel's Question 0.8.2 \cite{KeelPol} asks: Given a nef line bundle $L$
on a projective scheme $X$ over $\Fb_p$ such that $L$-equivalence
is bounded, is $L$ semi-ample? We now check that
on a smooth projective surface over any field, $L$-equivalence
is always bounded (Theorem \ref{Leq}). As a result, the examples in this paper
of nef but not semi-ample line bundles on smooth projective surfaces
over $\Fb_p$ give a negative answer to Question 0.8.2.

\begin{theorem}
\label{Leq}
For any nef line bundle $L$ on a normal proper algebraic space $X$ 
of dimension 2 over a field,
$L$-equivalence is bounded.
\end{theorem}

{\bf Proof. }First suppose $X$ is smooth and projective.
The theorem is clear if $L$ is numerically trivial, since
any two points lie on a curve. Also, if $L$ is big ($L^2>0$),
then there are only finitely many $L$-exceptional curves by
Lemma \ref{bound}.

So we can assume that $L^2=0$ and $L$ is numerically nontrivial.
By Lemma \ref{bound}, there are only finitely many $L$-exceptional
curves $A$ which are not numerically equivalent to a multiple of
$L$. There may be infinitely many curves $C$ which are numerically
equivalent to a multiple of $L$, but every such curve
is clearly disjoint from all other $L$-exceptional curves. Thus $L$-equivalence
is bounded.

Now let $X$ be normal, not necessarily smooth. There is a resolution
of singularities $f:M\arrow X$ with $M$ projective.
The $f^*(L)$-exceptional curves
on $M$ are the proper transforms of the $L$-exceptional curves together
with the finitely many curves contracted by $f$. This makes it clear
that $L$-equivalence is bounded on $X$ in the cases where $L$ is numerically
trivial or big. For $L^2=0$ and $L$ numerically nontrivial,
the previous paragraph's results applied to $f^*(L)$ on $M$
imply that
all but finitely many $L$-exceptional curves of $X$
are disjoint from the singular points
of $X$ and from all other $L$-exceptional curves; so again $L$-equivalence
is bounded. \qed

\section{General remarks on moving divisors}
\label{preparation}

Here we discuss the obstructions to moving codimension-one subvarieties
in general terms, to motivate our main technical result,
Theorem \ref{mainm}. The main application in this paper will
be a negative solution to Keel's Question 0.8.2 \cite{KeelPol}.

We want to find a curve $C$ with self-intersection zero
in some surface $M$ over $\Fb_p$ such that no multiple
of $C$ moves on $M$. In order for this to happen, $C$ must have
arithmetic genus at least 2
by Theorem \ref{genus1}.
Moreover, the normal bundle
of $C$ must be nontrivial by Lemma \ref{mumford}, stated by Mumford
and discussed below.
On the other hand,
the normal bundle of $C$ must be torsion in the Picard group of $C$,
because it is a line bundle of degree zero over $\Fb_p$
(Lemma \ref{torsion}). We might therefore
look for counterexamples in what seems to be the simplest remaining case:
where $C$ is a smooth curve of genus 2 whose normal bundle has order 2
in the Picard group of $C$. There are indeed
counterexamples of this type over $\Fb_p$ for any prime number $p$,
but it turns out to be simpler to give counterexamples where
$C$ is a curve of genus 2 with normal bundle of order $p$
equal to the characteristic. We give only the latter type
of counterexample in this paper.

The examples we give are curves of genus 2 in rational surfaces $M$,
obtained by blowing up $\P^1\times \P^1$
at 12 points on a smooth curve of genus 2
(so that the proper transform $C$ has $C^2=0$).
It is easy to choose the points we blow up to ensure
that the normal bundle of $C$ has order $p$ in the Picard group of $C$.
Theorem \ref{mainm} gives a Zariski open condition (on the family
of such arrangements of points and curves)
which implies that no multiple of $C$ moves on $M$.
That is the key point; a priori, to prove that no multiple
of $C$ moves would require us to check infinitely many open conditions.
In Theorem \ref{mainex}, we show that there are indeed
counterexamples to Keel's question of this type, over $\Fb_p$
for every prime number $p$.

Once we have a counterexample to Keel's question on a rational surface $M$,
we get counterexamples on many other surfaces, including
surfaces of general type. Take any smooth projective surface
$S$ with a surjective morphism $f:S\arrow M$; for example, take a ramified
covering of $M$ and resolve singularities. Let $L$ be a nef line bundle
on $M$ which is numerically nontrivial but has Iitaka dimension
$\kappa(M,L)\leq 0$, as in our examples. Then the pullback line bundle
$f^*L$ has the same properties. Indeed, we have 
$\kappa(S,f^*L)=\kappa(M,L)$ for every surjective morphism
$f:S\arrow M$ of normal varieties, by Ueno \cite[Theorem 5.13]{Ueno}.
Thus the line bundle $f^*L$ is a counterexample to Keel's
question on the surface $S$.

Here is the lemma stated by Mumford, as mentioned above. He left the proof
as a ``curiosity for the
reader'' \cite[Proposition, p.~336]{MumfordI}.
A proof was given by Ma\c sek \cite[Lemma, p.~682]{Masek}.
The lemma explains
our use of curves with nontrivial but torsion normal bundle
in order to give
counterexamples to Keel's question, but it will not actually
be used in the rest of the paper. We will state and prove
the lemma in any dimension, using the idea of killing
cohomology as in the proof of Theorem \ref{genus1}.

\begin{lemma}
\label{mumford}
Let $D$ be an effective Cartier divisor
in a normal projective variety $X$ over a field
of characteristic $p>0$,
and let $L$ be the line bundle $O(D)$ on $X$. 
Assume that the restriction of $L$ to $D$, the normal
bundle of $D$ in $X$, has finite order $m$ in the Picard
group of $C$ (so in particular the normal bundle is numerically
trivial). Then $L$ is semi-ample on $X$
(so that some positive
multiple of $D$ moves in its linear system with no base points)
if and only if the line bundle $mp^rL$ is trivial on 
the subscheme $mp^rD$ of $X$ for some $r\geq 0$.
\end{lemma}

Here, for any effective Cartier divisor $D$
in a normal scheme $X$, $D$ will be
defined locally by one equation $f=0$; then $aD$ denotes
the subscheme of $X$ defined locally by $f^a=0$, for any
positive integer $a$.
For any integer $a$, we will write $aL$ for the line bundle
$L^{\otimes a}=O(aD)$ on $M$.

In characteristic zero, it is not true that the triviality of the
line bundle $mL$ on the subscheme $mD$
implies that some multiple of $D$ moves in its linear system.
A counterexample is provided by the
$\P^1$-bundle over an elliptic curve associated to a nontrivial
extension of the trivial line bundle by itself. Here
a section $D$ is an elliptic curve with trivial normal bundle
(that is, $L:=O(D)$ is trivial on $D$) and yet no multiple of $D$ moves.
But there is an analogous result in characteristic zero, by Neeman.
Let $D$ be a curve in a smooth projective surface $X$ over a field $k$
of characteristic zero such that the normal bundle has finite order $m$.
Then some positive multiple of $D$ moves in its linear system on $X$
(or, equivalently, $L=O(D)$ is semi-ample on $X$) if and only if
the line bundle $mL$ is trivial on the subscheme $(m+1)D$
\cite[Article 2, Theorem 5.1]{Neeman}. Neeman works over the complex
numbers and considers smooth curves, but only minor changes to the argument
are needed to avoid those restrictions.

{\bf Proof of Lemma \ref{mumford}. }First suppose that
$L$ is semi-ample on $X$. Then $aL$ is basepoint-free on $X$ for
some positive integer $a$, and hence $aL$ is basepoint-free
on the subscheme $aD$. But $L$ is also numerically trivial
on $aD$, and so $aL$ is trivial on $aD$.
Clearly $a$ must be a multiple of $m$,
and so we can write $a=mp^rj$ for some $r\geq 0$ and some positive
integer $j$ which is not a multiple of $p$.
The group $\ker(\Pic(mp^rjD)\arrow \Pic(D))$ is $p$-primary by
the exact sequence used in the proof of Lemma \ref{mj} below.
Since the line bundle $mp^rL$ is trivial
on $D$ and $mp^rjL$ is trivial on $mp^rjD$, it follows
that $mp^rL$ is trivial on $mp^rjD$. A fortiori, $mp^rL$ is trivial
on $mp^rD$.

Conversely, suppose that $aL$ is trivial on $aD$
for $a=mp^r$ and some $r\geq 0$. By Lemma \ref{good},
to show that $L$ is semi-ample on $X$, it suffices to show
that $\kappa(X,L)$ is at least 1, which holds if
$H^0(X,nL)$ has dimension at least 2 for some positive integer $n$.
By the exact sequence
$$0\arrow O_X\arrow O_X(aD)\arrow O_X(aD)|_{aD}\arrow 0$$
of sheaves on $X$, we have an exact sequence of cohomology groups
$$H^0(X,O)\arrow H^0(X,aL)\arrow H^0(aD,aL)\arrow H^1(X,O).$$
Since $aL$ is trivial on $aD$, we have $h^0(X,aL)\geq 2$
and hence $L$ is semi-ample if $H^1(X,O)=0$. In general,
we use the idea of killing cohomology as in the proof
of Theorem \ref{genus1}. We find that there is a normal variety $Y$
with a finite morphism $f:Y\arrow X$ such that the obstruction
class in $H^1(X,O)$ pulls back to 0 in $H^1(Y,O)$. Therefore
$\kappa(Y,f^*L)\geq 1$. By Ueno \cite[Theorem 5.13]{Ueno},
it follows that $\kappa(X,L)\geq 1$. Thus $L$ is semi-ample.
\qed

\section{Obstruction theory for moving divisors}

In this section, we prove our main general result of obstruction
theory for moving divisors, Theorem \ref{mainm}.
It gives a sufficient condition on a divisor $D$ with torsion
normal bundle in a variety $M$ so that no multiple of $D$ moves.
(If the normal bundle is numerically trivial but not torsion,
then it is clear that no multiple of $D$ moves.)
The condition is Zariski open on the family of divisors with
normal bundle of a given order. 
Section \ref{orderpex} will apply this theorem to give counterexamples to 
Keel's Question 0.8.2.

\begin{theorem}
\label{mainm}
Let $D$ be a connected reduced
Cartier divisor in a normal projective variety $M$
over an algebraically closed field of characteristic $p>0$. Let $L$
denote the line bundle $O(D)$ on $M$. Suppose that
the restriction of $L$ to $D$ has order equal to $p$
in the Picard group of $D$ (so in particular $D$ has numerically
trivial normal bundle).
Suppose that the line bundle $pL$ is nontrivial
on the subscheme $2D$ of $M$. Finally, suppose that
the Frobenius maps
$F^*: H^1(D,-L)\arrow H^1(D,-pL)\cong H^1(D,O)$ and
$F^*:H^1(D,O)\arrow H^1(D,O)$ are injective.

Then $L$ is not semi-ample on $M$. More precisely,
no multiple of $D$ moves in its linear system on $M$.
\end{theorem}

{\bf Proof. }In view of the inclusions $H^0(X,O(D))\subset H^0(X,O(2D))
\subset \cdots$, it suffices to show that $H^0(X,O(p^{r+1}D))$ has dimension
1 for all $r\geq 0$. 
By the exact sequence of sheaves
$$0\arrow O_X\arrow O_X(p^{r+1}D)\arrow O_X(p^{r+1}D)|_{p^{r+1}D}\arrow 0,$$
it suffices to show that $H^0(p^{r+1}D,p^{r+1}L)=0$ for all
$r\geq 0$. 
That will be accomplished by Lemma \ref{mj}. 
The proof of Lemma \ref{mj} will require a simultaneous induction
to show that every regular function on $p^{r+1}D$
is constant and that $p^{r+1}L$ is nontrivial on $p^{r+1}D$.
\qed

\begin{lemma}
\label{mj}
We retain the assumptions of Theorem \ref{mainm}.
For any integer $j$ and any $r\geq 0$,
the group $H^0(p^rD,-p^rjL)$ is 0 for $j$ not a multiple
of $p$, and it injects into $H^0(D,-p^rjL)\cong k$ for $j$
a multiple of $p$.
Also, the line bundle $p^{r+1}L$ is nontrivial on $(p^r+1)D$
(and hence on $p^{r+1}D$).
\end{lemma}

{\bf Proof. }Since $D$ is connected and reduced,
the group $H^0(D,O)$ of regular functions on $D$ is equal to the
algebraically closed base field $k$.
For $r=0$, the first statement of the lemma holds because
$L$ is a line bundle of order equal to $p$ in the Picard
group of the divisor $D$. Also, $pL$ is nontrivial on $(p^0+1)D=2D$
by assumption.

For any positive integer $a$, we have an exact sequence of sheaves
supported on $D$,
$$0\arrow L^{\otimes -a}|_D \arrow O^*_{(a+1)D}\arrow O^*_{aD}\arrow 0.$$
Part of the long exact sequence on cohomology looks like:
$$H^1(D,-aL)\arrow \Pic((a+1)D)\arrow \Pic(aD).$$
Since the line bundle $pL$ is trivial on $D$, this sequence
shows that the isomorphism class of $pL$
on $2D$ is the image of some element $\eta$ in
$H^1(D,-L)$. Since $pL$ is nontrivial on $2D$, $\eta$ is not zero.

Now let $r\geq 1$ and assume the lemma for $r-1$. For any positive
integers $a$ and $b$, we have an exact sequence of sheaves supported on $D$,
$$0\arrow L^{\otimes -a}|_{bD}\arrow O_{(a+b)D}\arrow O_{aD}\arrow 0.$$
Tensoring with a line bundle and taking cohomology gives
the following exact sequence, for any
$1\leq i\leq p-1$ and any integer $j$:
$$0\arrow H^0(p^{r-1}D,-p^{r-1}(pj+i)L)\arrow H^0(p^{r-1}(i+1)D,
-p^{r}jL)\arrow H^0(p^{r-1}iD,-p^{r}jL).$$
Since $pj+i$ is not a multiple of $p$, the first $H^0$ group
is zero by induction. Combining these exact sequences
for $i=1,\ldots,p-1$, we find that $H^0(p^rD,-p^{r}jL)$ injects
into $H^0(p^{r-1}D,-p^rjL)$. 
By our inductive assumption again, the latter group
in turn restricts injectively to $H^0(D,-p^{r}jL)\cong k$.
This is the conclusion we want when $j$ is a multiple of $p$.

For $j$ not a multiple of $p$, we use that $p^rL$ is nontrivial
on $(p^{r-1}+1)D$, hence on $p^rD$. The exact sequence
$$H^1(D,-aL)\arrow \Pic((a+1)D)\arrow \Pic(aD)$$
for positive integers $a$, where $H^1(D,-aL)$ is a $k$-vector
space, shows that $\ker(\Pic(p^rD)
\arrow \Pic(D))$ is a $p$-primary group.
Therefore $p^rjL$ is nontrivial on $p^rD$ for all $j$
not a multiple of $p$. 
Since $H^0(p^rD,-p^rjL)$ injects
into $H^0(D,-p^rjL)\cong H^0(D,O)\cong k$, it follows that
$H^0(p^rD,-p^rjL)=0$ (since a nonzero section would
give a trivialization of $-p^rjL$ over $p^rD$). This is
what we want for $j$ not a multiple of $p$.

Recall that $\eta$ is the nonzero element of $H^1(D,-L)$
that describes the isomorphism class of the line bundle $pL$
on $2D$. By our assumptions on injectivity of Frobenius maps,
$(F^*)^{r}(\eta)$ is nonzero in $H^1(D,-p^rL)$. Consider the 
exact sequence
$$H^0((p^{r}+1)D,O^*)\arrow H^0(p^{r}D,O^*)\arrow H^1(D,-p^{r}L)
\arrow \Pic((p^{r}+1)D).$$
We have shown (by the case $j=0$, above) that the group
$H^0(p^rD,O)$ of regular functions on $p^rD$ consists
only of the constants $k$, and so the first arrow in this sequence
is surjective. Therefore the last arrow is injective, and so the
image of $(F^*)^r(\eta)$ in $\Pic((p^r+1)D)$ is nonzero. This means
that $p^{r+1}L$ is nontrivial on $(p^r+1)D$, completing the induction.
\qed

\section{Curves with normal bundle of order $p$}
\label{orderpex}

\begin{theorem}
\label{mainex}
For every prime number $p$, there is a smooth projective
surface $M$ over $\Fb_p$ and a smooth curve $C$ of genus 2
in $M$ such that $C$ has self-intersection zero and
no positive multiple of $C$ moves in its linear system
on $M$. Therefore the line bundle $L=O(C)$ is nef but
not semi-ample on $M$.
\end{theorem}

That is, Keel's Question 0.8.2 \cite{KeelPol} has a negative answer
in every characteristic $p>0$. The example here can be considered optimal,
in view of the positive results listed in section \ref{preparation}.

The surface $M$ will be rational. We will start
with a smooth curve of genus 2 in $X=\P^1\times\P^1$ and then blow up
to make its self-intersection zero. (It seems easier
to start with $\P^1\times\P^1$ rather than $\P^2$, since a curve
of geometric genus 2 in $\P^2$ is always singular.) The following
lemma is a first step. We say that a possibly singular curve $C$
over a field
of positive characteristic
is {\it ordinary }if the Frobenius map $F^*:H^1(C,O)\arrow H^1(C,O)$
is injective.

\begin{lemma}
\label{orderplemma}
Let $C$ be a smooth curve of bidegree $(2,3)$ in $\P^1\times \P^1$
(so $C$ has genus 2) over an algebraically closed field 
of characteristic $p>0$. Suppose that $C$ is ordinary. Then
there is an effective divisor $D$
of degree 12 on $C$ with the following property. Let $M$ be the blow-up
of $\P^1\times \P^1$ at the 12 points of $D$, and let $L$ be
the line bundle on $M$ associated to the proper transform of $C$.
Then $L$ has order $p$ on $C$, while
the line bundle $O(pC)=pL$
is nontrivial on the subscheme $2C$ of $M$.
\end{lemma}

{\bf Proof. }The tangent bundle of $X=\P^1\times \P^1$ restricted
to $C$ is an extension of two line bundles,
$0\arrow TC\arrow TX|_C\arrow N_{C/X}\arrow 0$. Let
$\beta\in H^1(C,TC-N_{C/X})$ be the class of this extension.
I claim that $\beta$ is highly nontrivial in the sense
that for each effective divisor
$B$ of degree 4 on $C$, the image of $\beta$
in $H^1(C,TC-N_{C/X}+B)$ is nonzero. 

We can view this image
as the class of the extension $0\arrow TC\arrow E\arrow N_{C/X}(-B)\arrow 0$,
where $E$ is a subsheaf of $TX|_C$ in an obvious way.
If this extension splits, then we have a nonzero map from the line
bundle $N_{C/X}(-B)$ to $E$ and hence to $TX|_C$. Here $N_{C/X}$
has degree 12 and so $N_{C/X}-B$ is a line bundle of degree 8.
On the other hand, since $C$ has bidegree $(2,3)$ and $T\P^1$
has degree 2, the restriction of $TX=T\P^1\oplus T\P^1$ to
$C$ is the direct sum of two line bundles of degrees 4 and 6.
So any map from a line bundle of degree 8 to $TX|_C$ is zero.
Thus the extension class
$\beta$ in $H^1(C,TC-N_{C/X})$ is highly nontrivial in the sense
claimed.

For the moment, fix a line bundle $L$ of order $p$ on $C$. (We know that
the $p$-torsion subgroup of $\Pic(C)$ is isomorphic to $(\Z/p)^2$,
because $C$ is ordinary of genus 2.) We want to blow up
$X=\P^1\times \P^1$ at a divisor $D$ of degree 12 on $C$ so that
the normal bundle of $C$ in the blow-up $M$ is isomorphic to $L$.
Here $L=N_{C/M}=N_{C/X}-D$, and so $D$ must be the divisor
of a nonzero section $\delta\in H^0(C,N_{C/X}-L)$. The line bundle
 $N_{C/X}-L$
has degree 12, and so this $H^0$ group has dimension 11 by Riemann-Roch;
thus there are many possible divisors $D$ such that $N_{C/M}$ is isomorphic
to $L$.

The restriction of the tangent bundle of the blow-up
$M$ to $C$ is the extension
$$0\arrow TC\arrow TM|_C\arrow N_{C/X}(-D)\arrow 0$$
which is the restriction of the extension $0\arrow TC
\arrow TX|_C \arrow N_{C/X}\arrow 0$ denoted $\beta$.
(That is, we can view $TM|_C$ as a subsheaf of $TX|_C$ using the
derivative of the map $M\arrow X$.) Equivalently, in terms
of an identification $L\cong N_{C/X}(-D)$, we can describe the
class of the extension $TM|_C$ in $H^1(C,TC-L)$ as the image
under the cup product
$$H^1(C,TC-N_{C/X})\otimes H^0(C,N_{C/X}-L)\arrow H^1(C,TC-L)$$
of $(\beta,\delta)$.

Knowing the class of $TM|_C$ as an extension (the Kodaira-Spencer
class) amounts to knowing
the isomorphism class of the subscheme $2C$ of $M$. Some references
are Morrow-Rossi \cite[Theorem 2.5]{MorR} or, in greater generality,
Illusie \cite[Theorem 1.5.1]{Illusieshort}. Therefore,
we can hope to describe the class of the line bundle $pL$
on $2C$ in terms of this extension, and this is accomplished
by Lemma \ref{pl}, a consequence of Illusie's results on deformation
theory. (Here we write $L$ for the line bundle
$O(C)$ on $M$. In fact, the line bundle $pL$ on $2C$ only depends on the
line bundle $L$ on $C$, because $pL$ is the Frobenius pullback
$F^*(L)$ via the map $F:2C\arrow C$.)
Since $L$ has order $p$
on $C$, the class of $pL$ on $2C$ lies in $\ker(\Pic(2C)\arrow
\Pic(C)) \cong H^1(C,-L)$, by the exact sequence in the proof
of Lemma \ref{mj}. We will use
Cartier's theorem that, for any smooth proper variety $X$ over an
algebraically closed field $k$ of characteristic $p>0$,
the $p$-torsion subgroup $\Pic(X)[p]$ tensored with $k$ injects naturally
into $H^0(X,\Omega^1)$ \cite[6.14.3]{Illusiecris}. 
Thus, for a curve $C$, $\Pic(C)[p]$
tensored with $k$ injects into $H^0(C,K_C)$.

We first recall Illusie's definition of the map $\Pic(X)[p]\arrow
H^0(X,\Omega^1)$ \cite[6.14]{Illusiecris}. We are given a line bundle $L$ 
with a trivialization of $L^{\otimes p}$. Let
$\{U_i\}$ be an open covering of $Y$ on which we choose
trivializations of $L$.
Let $e_{ij}\in O(U_i\cap U_j)^*$ be the transition functions;
then the trivialization of $L^{\otimes p}$ gives functions
$u_i\in O(U_i)^*$ such that $e_{ij}^p=u_j/u_i$. It follows
that $d\log u_j-d\log u_i=p\, d\log e_{ij}=0$, and so the 1-forms
$d\log u_i$ fit together to give the desired class
in $H^0(Y,\Omega^1)$.

The following lemma separates the two roles of $L$ in our problem:
a line bundle of order $p$ on $C$ and the normal bundle of $C$ in $M$.
This yields a clearer and more general statement.

\begin{lemma}
\label{pl}
Let $C$ be a smooth compact subvariety of a smooth variety $M$ over an
algebraically closed field
of characteristic $p>0$. Let $L$ be a line bundle on $M$
(or just on $2C$) such that $pL$ is trivial on $C$. Then
the class of the line bundle $pL$ in $\ker(\Pic(2C)\arrow \Pic(C))
\cong H^1(C,N_{C/M}^*)$ is the cup product of
the class of $L$ in $H^0(C,\Omega^1)$ with the class of the extension $TM|_C$
in $H^1(C,TC\otimes N_{C/M}^*)$:
$$H^0(C,\Omega^1_C)\otimes H^1(C,TC\otimes N_{C/M}^*)\arrow
H^1(C,N_{C/M}^*).$$
\end{lemma}

{\bf Proof. }This can be proved by an explicit cocycle calculation,
but we will instead deduce it from Illusie's general results.
Namely, let $G$ be a flat group scheme
over a scheme $S$. For any extension
$Y\inj Y'$ of a scheme $Y$ over $S$ by a square-zero ideal sheaf,
the obstruction to extending a $G$-torsor over $Y$ (in the fpqc
topology) to $Y'$
is the product of the Atiyah class of the $G$-torsor over $Y$
with the Kodaira-Spencer class of the extension $Y\inj Y'$
\cite[2.7.2]{Illusieshort}.

Apply this to the group scheme $G=\mu_p$ of $p$th roots
of unity over a field $k$ of characteristic $p>0$. Using the
exact sequence 
$$\begin{CD} H^1(Y,\mu_p)@>>> H^1(Y,G_m)@>>p>
H^1(Y,G_m),\end{CD}$$ 
we can rephrase the lemma in terms of the obstruction
to extending a $\mu_p$-torsor from $Y=C$ to $Y'=2C$.
The analogue
of the dual of the Lie algebra of $\mu_p$ in Illusie's theory
is the object $\mathfrak{g}^*=k[1]\oplus k$ of the derived category of $k$
\cite[Example 4.3.4]{Illusieshort}. This calculation uses that
$\mu_p$ is a codimension-one subgroup of the one-dimensional
group scheme $G_m$. 
As a result, the Atiyah class of a $\mu_p$-torsor on a smooth
scheme $Y$ over $k$ lies in
$\Ext^1_Y(\mathfrak{g}^*_Y,\Omega^1_Y)=H^0(Y,\Omega^1)\oplus H^1(Y,\Omega^1)$. The Kodaira-Spencer
class of the extension $Y\inj Y'$ lies in $H^1(Y,TY\otimes N_{Y/Y'}^*)$.
Finally, the obstruction to extending
a $\mu_p$-torsor lies in $H^1(Y,N_{Y/Y'}^*)\oplus H^2(Y,N_{Y/Y'}^*)$.
The second part of this obstruction is the obstruction to extending
the line bundle $L$ from $C$ to $2C$, which is zero since we are
given a line bundle $L$ on $2C$. The first part of this obstruction
is computed by Illusie's product formula. \qed (Lemma \ref{pl}).

We now return to the blow-up $M$ of $X=\P^1\times \P^1$. Let
$\beta$ in $H^1(C,TC-N_{C/X})$ be the class of the extension
$TX|_C$, and let $\delta$ in $H^0(C,N_{C/X}-L)$ be a section
whose zero set is the divisor on $C$ where we blow up.
Let $\gamma$
be the class of the $p$-torsion line bundle $L$ in $H^0(C,K_C)$.
By Lemma \ref{pl}
together with our earlier results, the class of $pL$ in
$\ker(\Pic(2C)\arrow \Pic(C))=H^1(C,-L)$ is the product $\beta\delta\gamma$ in:
$$H^1(C,TC-N_{C/X})\otimes H^0(C,N_{C/X}-L)\otimes H^0(C,K_C)
\arrow H^1(C,-L).$$
The lemma is proved if there is a line bundle $L$ on $C$ of order $p$ (which
determines $\gamma$) and a section
$\delta\in H^0(C,N_{C/X}-L)$ such that
the product $\beta\delta\gamma$ in $H^1(C,-L)$ is not zero. By Serre duality,
this product becomes:
$$\begin{CD} H^0(C,N_{C/X}-L)\otimes H^0(C,K_C)\otimes H^0(C,K_C+L)
\arrow H^0(C,2K_C+N_{C/X})
@>{\beta}>> k.$$
\end{CD}$$
We want to show: (*) for some line bundle $L$ on $C$ of order $p$ (which
determines $\gamma$), some $\delta\in H^0(C,N_{C/X}-L)$, and some
$\alpha\in H^0(C,K_C+L)$, the product $\beta(\delta\gamma\alpha)\in k$
is not zero.

Here $L$ is a nontrivial line bundle of degree 0 on the curve $C$ of genus 2,
and so $h^0(C,K_C+L)=1$ by Serre duality. Let $B_L$ be the base locus of
$K_C+L$ on $C$, that is, the zero locus of a nonzero section of $K_C+L$;
clearly $B_L$ is an effective divisor of degree 2 on $C$,
since $K_C+L$ has degree 2.
Then it is clear that the image of the product
$$H^0(C,N_{C/X}-L)\otimes H^0(C,K_C+L)\arrow H^0(C,K_C+N_{C/X})$$
is the subspace $H^0(C,K_C+N_{C/X}-B_L)\subset H^0(C,K_C+N_{C/X})$
of sections that vanish on $B_L$.

We can now use the assumption that
$C$ is ordinary to deduce that the finite group $\Pic(C)[p]\cong
(\Z/p)^2$ spans the $k$-vector space $H^0(C,K_C)\cong k^2$. Let $L_1$
and $L_2$ be two line bundles of order $p$ on $C$ whose classes
span $H^0(C,K_C)$. Let $B_1$ and $B_2$ be the base loci of $K_X+L_1$
and $K_X+L_2$, respectively, which are effective divisors of degree
2 on $C$. Then $B:=B_1+B_2$ is an effective divisor of degree 4
on $C$. For $i=1$ or 2, the image of the product
$$H^0(C,N_{C/X}-L_i)\otimes H^0(C,K_C+L_i)\arrow H^0(C,K_C+N_{C/X})$$
contains the subspace $H^0(C,K_C+N_{C/X}-B)\subset H^0(C,K_C+N_{C/X})$
of sections that vanish on $B$, by the previous paragraph.

As a result, the lemma is proved
if the product 
$$\begin{CD}
H^0(C,K_C+N_{C/X}-B)\otimes H^0(C,K_C)\arrow H^0(C,2K_C+N_{C/X}-B)
\subset H^0(C,2K_C+N_{C/X}) @>{\beta}>> k \end{CD}$$
is nonzero. Indeed, if that holds, then at least one of the two
order-$p$ line bundles $L_1$ or $L_2$ will have class in
$H^0(C,K_C)$ whose product with $H^0(C,K_C+N_{C/X}-B)$
has nonzero image under $\beta$, which implies the
statement (*) above.

By Castelnuovo's theorem \cite[p.~151]{ACGH}, the product map
$$H^0(C,K_C+N_{C/X}-B)\otimes H^0(C,K_C)\arrow H^0(C,2K_C
+N_{C/X}-B)$$
is surjective, because $C$ has genus $g$ at least 2 and
$K_C+N_{C/X}-B$ has degree at least $2g+1$ (in fact, it has degree
10, which is at least $2g+1=5$). So it remains only to show that
$\beta$ is nonzero on the subspace $H^0(C,2K_C+N_{C/X}-B)$ of
$H^0(C,2K_C+N_{C/X})$ for every effective divisor $B$ of degree 4.
This follows by Serre duality from the property of $\beta$
in $H^1(C,TC-N_{C/X})$ proved at the beginning of this proof:
the image of $\beta$ in $H^1(C,TC-N_{C/X}+B)$ is nonzero for every effective
divisor $B$ of degree 4 on $C$. \qed
(Lemma \ref{orderplemma}).

\begin{lemma}
\label{injective}
Let $(C,L)$ be a general pair over $\Fb_p$ with $C$ a smooth ordinary curve
of genus 2 and $L$ a line bundle of order $p$ on $C$. Then the Frobenius
map $F^*:H^1(C,L)\arrow H^1(C,pL)\cong H^1(C,O)$ is injective.
(Here $H^1(C,L)$ has dimension 1 and $H^1(C,O)$ has dimension 2.)
\end{lemma}

Note that the moduli space of pairs $(C,L)$ over $\Fb_p$
with $C$ ordinary of genus 2
and $L$ of order $p$ is irreducible, by the irreducibility
of the moduli space of curves of genus 2 together with
the theorem
that the geometric monodromy homomorphism for ordinary curves
(on the group
$(\Z/p)^2$ of line bundles killed by $p$) maps onto $GL_2(\Z/p)$
(Faltings-Chai \cite[Prop.~V.7.1]{FC}, Ekedahl \cite{Ekedahl}).
As a result, it makes sense to talk about
a general pair $(C,L)$, meaning the pairs outside some proper closed subset
of the moduli space.

{\bf Proof. }Let $C_0$ be
the union of two ordinary elliptic curves $E_1$ and $E_2$ identified at the
origin. Since $E_1$ is ordinary, there is a line bundle $L$ on $C_0$ 
that has order $p$
on $E_1$ and is trivial on $E_2$. Then $H^1(C_0,O)\cong H^1(E_1,O)
\oplus H^1(E_2,O)\cong k^2$. Also, $H^1(C_0,L)\cong H^1(E_1,L)
\oplus H^1(E_2,O)\cong k$. So the Frobenius map $F^*:H^1(C_0,L)
\arrow H^1(C_0,pL)\cong H^1(C,O)$ is the Frobenius map $H^1(E_2,O)\arrow
H^1(E_2,O)\subset H^1(E_1,O)\oplus H^1(E_2,O)$. Since $E_2$ is ordinary,
this Frobenius map is nonzero, hence injective.

We can deform $C_0$ to smooth ordinary curves of genus 2. In this
deformation, the Jacobians form a smooth family of ordinary
abelian surfaces. Therefore the line
bundle $L$ of order $p$ on $C_0$ can be deformed to a line bundle of order $p$
over the smooth curves (over an etale open subset of the parameter space).
It follows that $F^*:H^1(C,L)\arrow H^1(C,pL)$ is injective for some
pair $(C,L)$ with $C$ smooth ordinary of genus 2 and $L$ of order $p$,
hence for general such pairs. \qed

We can now give the first examples of nef but not semi-ample line
bundles on smooth projective varieties over $\Fb_p$.
We know that a general smooth curve $C$ of genus 2 is ordinary.
It follows from Lemma \ref{injective} that,
for a general ordinary smooth curve $C$ of genus 2,
every line bundle $L$ of order $p$ on $C$ has injective Frobenius map
$F^*: H^1(C,L)\arrow H^1(C,pL)\cong H^1(C,O)$. Let $C$ be an ordinary
smooth curve of genus 2 over $\Fb_p$ that has this injectivity property.

We can imbed this curve, like any smooth curve $C$ of genus 2,
as a curve of bidegree $(2,3)$
in $\P^1 \times \P^1$. Explicitly: use the line bundles $K_C$ and any line
bundle $A$ of degree 3 that is not of the form $K_C+p$ for any point $p$
in $C$. The line bundles $K_C$ and $A$ 
are both basepoint-free and have $h^0=2$,
and so they give two morphisms $C\arrow \P^1$,
hence a morphism $C\arrow \P^1\times \P^1$
of bidegree $(2,3)$. It is straightforward to check from the assumption
on $A$ that $C\arrow \P^1\times \P^1$ is an embedding.

By Lemma \ref{orderplemma}, there is a divisor $D$ of degree 12 on
$C$ over $\Fb_p$ 
such that blowing up $\P^1\times \P^1$ at the divisor $D$ gives
a surface $M$ with the following property. The line bundle $L:=O(C)$ on $M$ 
has order $p$ when restricted to $C$, while $pL$ is nontrivial on $2C$.
(Here $C$ denotes the proper transform in $M$ of the curve with the same
name in $\P^1\times \P^1$.)

We have arranged for all the hypotheses of Theorem \ref{mainm}.
Therefore $L=O(C)$ is not semi-ample on $M$, and in fact
no multiple of $C$ moves on $M$. \qed

\section{Nef and big but not semi-ample}

To conclude, it is easy to use our examples in dimension 2 to produce
a similar example in dimension 3, but now involving a nef
{\it and big }line bundle. In the simplest example, the 3-fold
is rational. By the argument in section \ref{preparation},
this gives examples on many other varieties, in particular
on 3-folds of general type.

\begin{theorem}
\label{bigex}
For any prime number $p$,
there is a nef and big line bundle $L$ on a smooth projective
3-fold $W$ over $\Fb_p$ which is not semi-ample.
\end{theorem}

Equivalently, by Theorem \ref{good}, the ring
$R(W,L)=\oplus_{a\geq 0} H^0(W,aL)$ is not finitely generated
over $\Fb_p$. There is no example as in Theorem \ref{bigex}
in dimension 2, by Artin \cite[proof of Theorem 2.9(B)]{Artin}
or Keel \cite{KeelAnn}.

{\bf Proof. }By Theorem \ref{mainex}, for every prime
number $p$, there is a smooth
projective surface $X$ over $k=\Fb_p$ and a nef line bundle $L_1$
on $X$ which is not semi-ample. Equivalently, as mentioned
in section \ref{notation}, if we let $L_2$ be an ample
line bundle on $X$, then the ring $R(X,L_1,L_2):=
\oplus_{a,b\geq 0}H^0(X,aL_1+bL_2)$ is not finitely
generated over $k$. Since $L_2$ and $L_1+L_2$ are ample,
this ring has Iitaka dimension 3, meaning that the subspaces
of total degree $d$ grow at least like a positive constant times $d^3$.

Let $W$ be the projective bundle $P(L_1\oplus L_2)$ of hyperplanes
in $L_1\oplus L_2$. The line bundle $O(1)$ on the $\P^1$-bundle
$\pi:W\arrow X$ is easily checked to be nef, since $L_1$ and $L_2$
are nef. We have $\pi_*O(1)=L_1\oplus L_2$ and more
generally $\pi_*O(d)=S^d(L_1\oplus L_2)=\oplus_{i=0}^d L_1^{\otimes i}
\oplus L_2^{\otimes d-i}$. Therefore
$$R(W,O(1))=R(X,L_1,L_2).$$
So the ring $R(W,O(1))$ has Iitaka dimension 3 but is not finitely
generated. Therefore the nef line bundle $O(1)$ on the 3-fold $W$
is big but not semi-ample. \qed


\small \sc DPMMS, Wilberforce Road,
Cambridge CB3 0WB, England

b.totaro@dpmms.cam.ac.uk
\end{document}